\documentclass[12pt]{article}
\usepackage{amssymb}
\usepackage{amsmath}
\usepackage{epsfig}
\usepackage{theorem}

\newcommand{\bA}{{\mathbb A}}
\newcommand{\bC}{{\mathbb C}}

\newcommand{\bG}{{\mathbb G}}
\newcommand{\bQ}{{\mathbb Q}}
\newcommand{\bR}{{\mathbb R}}
\newcommand{\bP}{{\mathbb P}}
\newcommand{\bZ}{{\mathbb Z}}

\newcommand{\cB}{{\mathcal B}}
\newcommand{\cC}{{\mathcal C}}
\newcommand{\cD}{{\mathcal D}}
\newcommand{\D}{{\mathcal D}}

\newcommand{\cG}{{\mathcal G}}

\newcommand{\cL}{{\mathcal L}}

\newcommand{\cO}{{\mathcal O}}

\newcommand{\cT}{{\mathcal T}}
\newcommand{\cU}{{\mathcal U}}
\newcommand{\cX}{{\mathcal X}}
\newcommand{\X}{{\mathcal X}}

\newcommand{\Z}{{\mathcal Z}}
\newcommand{\fp}{{\mathfrak p}}
\newcommand{\fq}{{\mathfrak q}}
\newcommand{\fo}{{\mathfrak o}}
\newcommand{\Bl}{\mathrm{Bl}}
\newcommand{\Br}{\mathrm{Br}}
\newcommand{\GL}{\mathrm{GL}}
\newcommand{\Id}{\mathrm{Id}}
\newcommand{\Tr}{{\mathrm{Tr}}}
\newcommand{\Spec}{\mathrm{Spec}}
\newcommand{\rank}{\mathrm{rank}}
\newcommand{\Proj}{\mathrm{Proj}}
\newcommand{\Pic}{\mathrm{Pic}}
\newcommand{\End}{\mathrm{End}}
\newcommand{\ra}{\rightarrow}

{\theorembodyfont{\rm}
\newtheorem{dfn}{Definition}[section]
\newtheorem{coro}[dfn]{Corollary}

\newtheorem{exam}[dfn]{Example}
\newtheorem{rem}[dfn]{Remark}

\newtheorem{prob}[dfn]{Problem}
}

\newtheorem{conj}[dfn]{Conjecture}
\newtheorem{theo}[dfn]{Theorem}
\newtheorem{prop}[dfn]{Proposition}

\author{Brendan Hassett and Yuri Tschinkel}
\title{\Large Density of integral points on algebraic varieties}

\date{\today}

\begin{document}

\maketitle

\section{Introduction}

Let $K$ be a number field, $S$ a finite
set of valuations of $K$, including
the archimedean valuations,
and  ${\cal O}_S$ the ring of $S$-integers.
Let $X$ be an algebraic variety defined over $K$ and
$D$ a divisor on $X$.  We will use
$\X$ and $\D$ to denote models
over $\Spec({\cal O}_S)$.

We will say that integral points on $(X,D)$
(see Section \ref{sect:generalities}
for a precise definition)
are potentially dense if they are Zariski dense
on some model $(\cX,\cD)$, after a finite
extension of the ground field and after enlarging $S$.
A central problem in arithmetic geometry is
to find conditions insuring potential density (or nondensity)
of integral points.
This question motivates many
interesting and concrete problems in
classical number theory, transcendence theory
and algebraic geometry, some of which will be presented below.

If we think about general reasons for the density of points - the
first idea would be to look for the
presence of a large automorphism group.
There are many beautiful examples both for rational and integral
points, like K3 surfaces
given by a bihomogeneous
$(2,2,2)$ form in $\bP^1\times \bP^1\times \bP^1 $
or the classical Markov equation $x^2 + y^2 +z^2 = 3xyz$.
However, large automorphism groups
are ``sporadic'' - they are hard to
find and usually, they are not well behaved in families.
There is one notable exception - namely automorphisms of
algebraic groups, like tori and abelian varieties.

Thus it is not a surprise that the
main geometric reason for the abundance of rational points
on varieties treated in the recent papers \cite{harris-tschi},
\cite{bog-tschi}, \cite{hassett-tschi}
is the presence of elliptic or, more generally, abelian
fibrations with {\em multisections}
having a dense set of rational points and
subject to some {\em nondegeneracy} conditions.
Most of the effort goes into ensuring
these conditions.

In this paper we focus on cases when $D$ is nonempty.
We give a systematic treatment of known approaches
to potential density
and present several new ideas for proofs.
The analogs of elliptic fibrations
in log geometry are conic bundles with a bisection removed.  We develop
the necessary techniques to translate the presence of
such structures to statements
about density of integral points and
give a number of applications.

The paper is organized as follows:
in Section \ref{sect:generalities}
we introduce the main definitions and notations.
Section \ref{sect:morphisms}
is geometrical - we introduce the relevant concepts from the
log minimal model program and formulate several geometric problems
inspired by questions about integral points.
In Section \ref{sect:fibration}, we recall the fibration method
and nondegeneracy properties of multisections.  We consider approximation
methods in Section \ref{sect:approx}.  Section \ref{sect:conic}
is devoted to the study of integral points on conic bundles with
sections and bisections removed.
In the final section, we survey the known results concerning
potential density for integral point on log K3 surfaces.

\

{\bf Acknowledgements.}
The first author was partially supported by an NSF
Postdoctoral Research Fellowship.  The
second author was partially supported by the NSA.
We benefitted from conversations 
with Y. Andr\'e, F. Bogomolov,
A. Chambert-Loir, J.-L. Colliot-Th\'el\`ene, J. Koll\'ar,
D. McKinnon, and B. Mazur.  We are
grateful to P. Vojta for comments that
improved the paper, especially Proposition \ref{prop:nocover},
and to D.W. Masser for information on specialization
of nondegenerate sections.  
Our approach in Section \ref{sect:conic} is inspired by the
work of F. Beukers (see \cite{Beukers1} and \cite{Beukers2}).

\section{Generalities}
\label{sect:generalities}

\subsection{Integral points}

Let $\pi \,:\, {\cal U}\ra \Spec({\cal O}_S)$
be a flat scheme
over ${\cal O}_S$ with generic fiber $U$.
An integral point on ${\cal U}$ is a section of $\pi$; the
set of such points is denoted ${\cal U}({\cal O}_S)$.

In the sequel, ${\cU}$ will be the
complement to a reduced effective
Weil divisor $\D$ in a normal proper scheme $\X$, both generally
flat over $\Spec({\cO}_S)$. 
Hence an $S$-integral point $P$ of $(\X,\D)$ is a section
$s_P\,:\, \Spec({\cO}_S)\ra \X$
of $\pi$, which does not intersect $\D$, that is,
for each prime ideal ${\fp}\in
\Spec({\cO}_S)$ we have $s_P({\fp})\notin
\D_{\fp}$.
We denote by
$X$ (resp. $D$) the corresponding generic fiber.
We generally assume that $X$ is a variety
(i.e., a geometrically integral scheme);
frequently $X$ is smooth and $D$ is normal crossings.
Potential density of integral points on $(\X,\D)$
does not depend on the choice of $S$ or
on the choices of models over $\Spec(\cO_S)$,
so we will not always specify them.
Hopefully, this will not create any confusion.

If $D$ is empty then every $K$-rational point of $X$ is
an $S$-integral point for $(\X,\D)$ (on some model).
Every $K$-rational point of $X$, not
contained in $D$ is $S$-integral on $(\X,\D)$
for $S$ large enough.
Clearly, for any $\X$ and $\D$ there exists a finite extension
$K'/K$ and a finite set $S'$ of prime ideals in ${\cO}_{K'}$
such that there is an $S'$-integral point on $(\X', \D')$
(where $\X'$ is the basechange of $\X$ to $\Spec({\cO}_S')$).

The definition of integral points can be generalized
as follows: let $\Z$ be any subscheme of $\X$, flat
over ${\cO}_S$.
An $S$-integral point for $(\X,\Z)$ is an $\cO_S$-valued point
of $\X\setminus \Z$.

\subsection{Vojta's conjecture}

A {\em pair} consists of a proper normal variety $X$ 
and a reduced effective Weil divisor $D\subset X$.  
A {\em morphism of pairs} $\varphi:(X_1,D_1)\ra (X_2,D_2)$ is
a morphism $\varphi:X_1\ra X_2$ such that $\varphi^{-1}(D_2)$
is a subset of $D_1$.  In particular, $\varphi$ restricts to
a morphism $X_1\setminus D_1 \ra X_2 \setminus D_2$.
A morphism of pairs is {\em dominant} 
if $\varphi:X_1\ra X_2$ is dominant.  If $(X_1,D_1)$
dominates $(X_2,D_2)$ then integral points are
dense on $(X_2,D_2)$ when they are dense on $(X_1,D_1)$
(after choosing appropriate integral models.)
A morphism of pairs is {\em proper} if $\varphi:X_1 \ra X_2$
is proper and the restriction $X_1 \setminus D_1 \ra X_2 \setminus D_2$
is also proper;  equivalently, we may assume that 
$\varphi:X_1 \ra X_2$ is proper and $D_1$ is a subset of $\varphi^{-1}(D_2)$.
A {\em resolution} of the pair $(X,D)$
is a proper morphism of pairs
$\rho:({\tilde X},{\tilde D}) \ra (X,D)$
such that $\rho:{\tilde X} \ra X$ is birational,
$\tilde X$ is smooth, and $\tilde D$ is normal crossings.

Let $X$ be a normal proper variety of dimension $d$.  Recall that
a Cartier divisor $D\subset X$ is {\em big} if 
$h^0(\cO_X(nD))>C n^d$
for some $C>0$ and all $n$ sufficiently large 
and divisible.  

\begin{dfn}
A pair $(X,D)$ is of log general type if it admits a
resolution
$\rho:({\tilde X},{\tilde D}) \ra (X,D)$
with $\omega_{\tilde{X}}(\tilde{D})$ big.
\end{dfn}
Let us remark that the definition does not depend on the
resolution.

\begin{conj}(Vojta, \cite{vojta})
Let $(X,D)$ be a pair of log general type.  Then
integral points on $(X,D)$ are not potentially dense.
\end{conj}
This conjecture
is known when $X$ is a semiabelian variety (\cite{faltings},
\cite{vojta-1}, \cite{mcquillan}).
Vojta's conjecture implies that a pair with dense integral points
cannot dominate a pair of log general type.

We are interested in geometric conditions which
would insure potential density of integral points.
The most naive statement would be the direct converse
to Vojta's conjecture. However this can't be true
even when $D=\emptyset$.
Indeed, varieties which are not of general type may
dominate varieties of general type, or more generally,
admit finite \'etale covers which dominate varieties of
general type (see the examples in \cite{colliot-swd}).
In the next section we will analyze other types of covers with
the same arithmetic property.

\section{Geometry}
\label{sect:morphisms}

\subsection{Morphisms of pairs}
\label{sect:morphisms-1}

\begin{dfn}
We will say that a class of
dominant morphisms of pairs
$
\varphi\,:\, (X_1,D_1)\ra (X_2,D_2)
$
is arithmetically continuous
if the density of integral points on
$(X_2,D_2)$ implies potential
density of integral points on
$(X_1,D_1)$.
\end{dfn}

For example, assume that $D=\emptyset$. Then any projective bundle
in the Zariski topology $\bP\ra X$ is arithmetically continuous.
In the following sections we present other examples of
arithmetically continuous morphisms of pairs.

\begin{dfn}
A pseudo-\'etale cover of pairs
$\varphi\,:\,
(X_1,D_1)\ra (X_2,D_2)$ is a 
proper dominant morphism of pairs such that

a) $\varphi:X_1 \ra X_2$ is generically finite, and

b) the map from the normalization $X_2^{\rm norm}$ of $X_2$ (in the function
field of $X_1$) onto $X_2$ is \'etale away from $D_2$.
\end{dfn}

\begin{rem}\label{rem:nc}
For every pair $(X,D)$ 
there exists a birational pseudo-\'etale morphism
$\varphi\,:\, (\tilde{X},\tilde{D})\ra (X,D)$
such that $\tilde{X}$ is smooth and $\tilde{D}$ is normal
crossings.
\end{rem}

The following theorem is a formal generalization of
the well-known theorem of Chevalley-Weil.
It shows that potential density is
stable under pseudo-\'etale covers of pairs.

\begin{theo}
\label{theo:chevalley-weil}
Let $\varphi\,:\,
(X_1,D_1)\ra (X_2, D_2)$ be a pseudo-\'etale cover
of pairs. Then $\varphi$ is arithmetically continuous.
\end{theo}

\begin{rem}
An elliptic fibration $E\ra X$, isotrivial on
$X\setminus D$, is
arithmetically continuous.
Indeed, it splits after a pseudo-\'etale
morphism of pairs and we can apply Theorem \ref{theo:chevalley-weil}.
\end{rem}

The following example is an integral analog of
the example of Skorobogatov, Colliot-Th\'el\`ene and
Swinnerton-Dyer (\cite{colliot-swd})
of a variety which does not dominate a variety
of general type but admits an
\'etale cover which does.

\begin{exam}
Consider $\bP^1\times \bP^1$ with coordinates $(x_1,y_1), (x_2,y_2)$
and involutions
$$j_1(x_1,y_1)=(-x_1,y_1) \quad j_2(x_2,y_2)=(y_2,x_2)$$
on the factors.
Let $j$ be the induced involution on the product;  it has
fixed points
$$
\begin{array}{cccr}
x_1=0  &  & x_2 = & y_2\\
x_1=0  &  &  x_2 = & -y_2 \\
y_1=0  &  & x_2 = & y_2 \\
y_1 =0 &  & x_2 = & -y_2  \end{array} .
$$
The first projection induces a map of quotients
$$({\bP}^1\times {\bP^1})/\left< j \right>
\ra {\bP}^1/\left<j_1 \right>.$$
We use $X$ to denote the source;  the target is 
just $\Proj({\bC}[x_1^2,y_1])\simeq {\bP}^1.$
Hence we obtain a fibration $f:X \ra {\bP}^1$.  
Note that $f$ has two nonreduced fibers, corresponding to
$x_1=0$ and $y_1=0$ respectively.  
Let $D$ be the image in $X$ of
$$
(x_1=0 )\cup (y_1 = 0) \cup (x_2=m_2y_2 ) \cup (x_2 =m_1y_2)
$$
(where $m_1,m_2$ are nonzero and distinct).
Since $D$ intersects the
general fiber of $f$ in just two points, $(X,D)$ is not of log general type.

We can represent $X$ as a degenerate quartic Del Pezzo surface
with four A1 singularities (see figure~\ref{logsurface1}).
\begin{figure}
\centerline{\psfig{figure=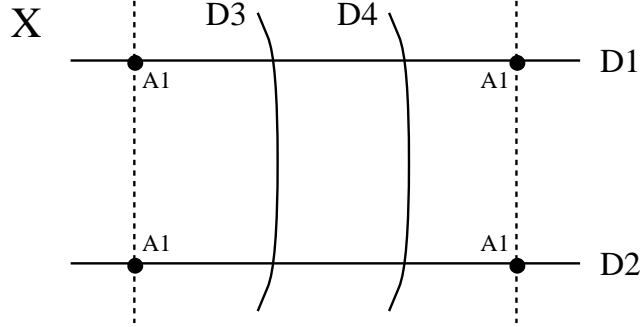}}
\caption{The log surface $(X,D)$}\label{logsurface1}
\end{figure}
If we fix invariants
$$
a=x_1^2x_2y_2, \ b=x_1^2(x_2^2+y_2^2), \ c=x_1y_1(x_2^2-y_2^2),
\ d=y_1^2(x_2^2+y_2^2),\ e=y_1^2x_2y_2
$$
then $X$ is given as a complete intersection of
two quadrics:
$$
ad=be, \hskip 1cm c^2=bd-4ae.
$$
The components of $D$ satisfy the equations
\begin{eqnarray*}
D_1&=&\{a=b=c=0\} \\
D_2&=& \{c=d=e=0\} \\
D_3&=& \{  (1+m_1^2)a-m_1b=(1+m_1^2)e-m_1d=0 \} \\
D_4&=& \{(1+m_2^2)a-m_2b=(1+m_2^2)e-m_2d=0 \}.
\end{eqnarray*}

We claim that
$(X,D)$ does not admit a dominant map onto
a variety of log general type and
that there exists a pseudo-\'etale
cover of $(X,D)$ which does.
Indeed, the preimage of $X\setminus D$ in
$\bP^1\times \bP^1$ is
$$
(\bA^1\setminus 0)\times (\bP^1 \setminus \{m_1,m_2, 1/m_1,
1/m_2\}),
$$
which dominates a curve of log general type,
namely, $\bP^1$ minus four
points.
However, $(X,D)$ itself cannot dominate a curve of
log general type.  Any such curve must be rational,
with at least three points removed;  however, the boundary $D$
contains at most two mutually disjoint irreducible components.
\end{exam}

The following was put forward as a possible converse to
Vojta's conjecture.

\begin{prob}[Strong converse to Vojta's conjecture]
\label{prob:main}
Assume that the pair $(X_2,D_2)$ does not admit
a pseudo-\'etale cover
$(X_1,D_1)\ra (X_2,D_2)$ such that
$(X_1,D_1)$ dominates a pair of log general type. Are
integral points for $(X_2,D_2)$ potentially dense?
\end{prob}

\subsection{Projective bundles in the \'etale topology}

We would like to produce further classes
of dominant arithmetically continuous morphisms
$(X_1,D_1) \ra (X_2,D_2)$.

\begin{theo}
Let $\varphi \, :\, (X_1, D_1 ) \ra (X_2,D_2)$ be a
projective morphism of pairs such
that $\varphi$ is a projective bundle (in the \'etale
topology) over $X_2  \setminus D_2$.
We also assume that $\varphi^{-1}(D_2) = D_1$.
Then $\varphi$ is arithmetically continuous.
\end{theo}
{\em Proof.}  We are very grateful to Prof. Colliot-Th\'el\`ene for
suggesting this proof.

Choose models $(\X_i,\D_i)$ ($i=1,2$) over
some ring of integers $\cO_S$, so that the morphism $\varphi$
is well-defined and satisfies our hypotheses.  (We enlarge $S$
as necessary.)

We recall basic properties of the Brauer group $\Br(\cO_S)$.
Let $v$ denote a place for the quotient field $K$
and $K_{v}$ the corresponding completion.  Classfield
theory gives the following exact sequence
$$0 \ra \Br(\cO_S)\ra \Br(K) \ra \oplus_{v \not \in S} \Br(K_{v}).$$
The Brauer groups of the local fields
corresponding to nonarchimedean valuations are isomorphic to $\bQ / \bZ$.
Given a finite extension of $K_{w}/K_{v}$ of
degree $n$ the induced map on Brauer groups is multiplication by
$n$.

Each $\cO_S$-integral point of $(X_2,D_2)$ yields an element
of  $\Br(\cO_S)$ of order $r$. This gives elements of $ \Br(K_{v})$
which are zero unless $v\in S$.
It suffices to find an extension $K'/K$ inducing a cyclic extension
of $K_v$ of order divisible by $r$ for all $v\in S$.
Indeed, such an extension necessarily kills any element of
$\Br(\cO_S)$ of order $r$.\hfill $\square$

\begin{rem}
Let $X$ be a smooth simply connected projective variety
which does not dominate a variety of general type.
It may admit an projective bundle (in the \'etale topology)
$\varphi\,:\, {\bP}\ra X$,
for example
if $X$ is a K3 surface.
However, $\bP$ cannot
dominate a variety of general type.
Indeed, given a dominant morphism $\pi\,:\, \bP\ra Y$,
the fibers of $\varphi$
are mapped to points by $\pi$.
In particular, $\pi$ necessarily factors through $\varphi$.
(We are grateful to J. Koll\'ar for emphasizing this point.)
\end{rem}

\begin{prob} [Geometric counterexamples
to Problem \ref{prob:main}]
\label{prob:eee}
Are there pairs which do not admit
pseudo-\'etale covers
dominating pairs of
log general type but which do admit arithmetically continuous
covers dominating pairs of log general type?
\end{prob}

\subsection{Punctured varieties}

In Section \ref{sect:morphisms-1} we have seen that
potential density of integral points is preserved under
pseudo-\'etale covers.
It is not an easy task, in general, to check whether or not some given
variety (like an elliptic surface) admits a (pseudo-) \'etale
cover dominating a variety of general type.
What happens if we modify the variety (or pair) without changing
the fundamental group?

\begin{prob}[Geometric puncturing problem]
\label{prob:geom}
Let $X$ be a projective variety with
canonical singularities and $Z$ a subvariety of codimension
$\ge 2$.
Assume that no (pseudo-) \'etale cover of $(X,\emptyset)$
dominates a variety  of general type.
Then $(X,Z)$ admits no pseudo-\'etale covers
dominating a pair of log general type.
A weaker version would
be to assume that $X$ and $Z$ are smooth.
\end{prob}
By definition, a pseudo-\'etale cover of $(X,Z)$
is a pseudo-\'etale cover of a pair $(X',D')$,
where $X'$ is proper over $X$ and
$X'\setminus D' \simeq X\setminus Z$.

\begin{prop}\label{prop:nocover}
Assume $X$ and $Z$ are as in Problem~\ref{prob:geom},
and that $X$ is smooth.
Then

a) No pseudo-\'etale covers of $(X,Z)$ dominate
a curve of log general type.

b) No pseudo-\'etale covers of $(X,Z)$ dominate
a variety of log general type of the same dimension.
\end{prop}

{\em Proof. } 
Suppose we have
a pseudo-\'etale cover $\rho : (X_1,D_1)\ra(X,Z)$
and a dominant morphism $\varphi : (X_1,D_1) \ra (X_2,D_2)$ to
a variety of log general type.  By Remark~\ref{rem:nc},
we may take the $X_i$ smooth and the $D_i$ normal crossings.
Since $D_1$ is exceptional with respect to $\rho$,
Iitaka's Covering Theorem (\cite{iitaka} Theorem 10.5)
yields an equality of Kodaira dimensions
$$
\kappa(K_X)=\kappa(K_{X_1}+D_1).
$$  

Assume first that $X_2$ is a curve.  We claim it has
genus zero or one.  Let $X^{\rm norm}$ be the normalization
of $X$ in the function field of $X_1$.  The induced morphism
$g:X^{\rm norm}\ra X$ is finite, surjective, and branched
only over $Z$, a codimension $\ge 2$ subset of $X$.  Since
$X$ is smooth, it follows that $g$ is \'etale.  
If $X_2$ has genus $\ge 2$ then $\varphi:X_1\ra X_2$ 
is constant along the fibers of $X_1\ra X^{\rm norm}$,
and thus descends to a map $X^{\rm norm} \ra X_2$.  
This would contradict our assumption that no \'etale cover of $X$
dominates a variety of general type.  

Choose a point $p\in D_2$ and 
consider the divisor $F=\varphi^{-1}p$.
Note that $2F$ moves because $2p$ moves on $X_2$.  
However, $2F$ is supported in $D_1$,
which lies in the exceptional locus for $\rho$,
and we obtain a contradiction.  

Now assume $\varphi$ is generically finite.  
We apply the Logarithmic Ramification Formula to
$\varphi$ (see \cite{iitaka} Theorem 11.5)
$$K_{X_1}+D_1=\varphi^*(K_{X_2}+D_2)+R$$
where $R$ is the (effective)
logarithmic ramification divisor.  
Applying the Covering Theorem again, we find
that $\kappa(K_{X_1}+D_1-R)=\kappa(K_{X_2}+D_2)=\dim(X)$.
It follows that $K_{X_1}+D_1$ is also big, which contradicts
the assumption that $X$ is not of general type.  
\hfill $\square$

\begin{prob}[Arithmetic puncturing problem]
\label{prob:punc}
Let $X$ be a projective variety with canonical singularities
and $Z$ a
subvariety of codimension $\ge 2$.  Assume that rational points
on $X$ are potentially dense.  Are integral points on
$(X,Z)$ potentially dense?
\end{prob}
For simplicity, one might first assume that $X$ and $Z$
are smooth.
Note that some conditions on the singularities of
$X$ are necessary. For example, blow up $\bP^2$ in
20 points lying along a smooth quartic curve $C$.
Assume that the divisor class of the points equals $5H$,
where $H$ is the hyperplane class of $C$.
Then the linear series of quintics containing the 20 points
gives a birational map
contracting $C$. Let $X$ be the resulting surface and $Z$ the
singular point. Rational points on $X$ are dense but
density of integral points
on $(X,Z)$ would contradict Vojta's conjecture.

\begin{rem}
Assume that Problem \ref{prob:punc}
has a positive solution.
Then potential density of
rational points holds for all K3 surfaces.

Indeed, if $Y$ is a K3 surface of degree $2n$ then
potential density of
rational points holds for the
symmetric product $X=Y^{(n)}$ (see \cite{hassett-tschi}).
Denote by $Z$ the large diagonal in $X$ and by
$\Delta$ the large diagonal in $Y^n$ (the ordinary product).
Assume
that integral points on $(X, Z)$ are potentially dense.
Then, by Theorem~\ref{theo:chevalley-weil}
integral points on $(Y^n, \Delta)$ are potentially dense.
This implies potential density for
rational points on $Y$.
\end{rem}

\section{The fibration method and nondegenerate multisections}
\label{sect:fibration}

This section is included as motivation.
Let $B$ be an algebraic variety, defined over
a number field $K$ and  $\pi \,:\, G\ra B$
be a group scheme over $B$.
We will be mostly
interested in the case when the generic fiber is
an abelian variety or a split torus $\bG_m^n$.
Let $s$ be a section  of $\pi$.
Shrinking the base we may assume that all fibers of $G$
are smooth.
We will say that
$s$ is {\em nondegenerate} if $\cup_n s^n$ is Zariski dense in
$G$.

\begin{prob}[Specialization]\label{prob:nondegenerate}
Assume that $G\ra B$ has a nondegenerate section $s$.
Describe the set of $b\in B(K)$ such that
$s(b)$ is nondegenerate in the fiber $G_b$.
\end{prob}

For simple abelian varieties over a field a point of infinite
order is nondegenerate.
If $E\ra B$ is a Jacobian elliptic fibration
with a section $s$ of infinite order then
this section is automatically nondegenerate, and
$s(b)$ is nondegenerate if it is nontorsion.
By a result of N\'eron
(see \cite{serre-MordellWeil} 11.1),
the set of $b\in B(K)$ such that
$s(b)$ is not of infinite order is {\em thin};  this
holds true for abelian fibrations of arbitrary dimension.  

For abelian fibrations $A\ra B$ with higher-dimensional fibers,
one must also understand how rings of endomorphisms
specialize.  
The set of $b\in B(K)$ for which the restriction 
$$\End(A) \ra \End(A(b))$$
fails to be surjective is also thin;  this is
a result of Noot \cite{noot} Corollary 1.5.
In particular, a nondegenerate section
of a family of generically simple abelian varieties
specializes to a nondegenerate point outside a thin set of fibers.

More generally, given an arbitrary abelian fibration $A\ra B$
and a nondegenerate section $s$, the set of $b\in B(K)$ such that
$s(b)$ is degenerate is thin in $B$.  (We are grateful to
Masser for pointing out the proof.)  After replacing $A$
by an isogenous abelian variety and taking a finite
extension of the function field $K(B)$, we obtain a
family $A'\ra B'$ with 
$A'\simeq A^{r_1}_1 \times \ldots \times A^{r_m}_m$, where
the $A_j$ are (geometrically) simple and mutually non-isogenous.
By the Theorems of N\'eron and Noot, the $A_j(b')$
are simple and mutually non-isogenous away from some thin
subset of $B'$.  A section $s'$ of $A'\ra B'$ is
nondegenerate iff its projection onto each factor $A^{r_j}_j$
is nondegenerate;  for $b'$ not contained in
our thin subset, $s'(b')$ is nondegenerate iff its
projection onto each $A^{r_j}_j(b')$ is nondegenerate.
Hence we are reduced to proving the claim for each
$A^{r_j}_j$.  Since $A_j$ is simple,
a section $s_j$ of $A^{r_j}_j$ is 
nondegenerate iff its projections $s_{j,1},\ldots,s_{j,r_j}$
are linearly independent over $\End(A_j)$.  Away from
a thin subset of $B'$, the same statement holds for
the specializations to $b'$.  However, N\'eron's theorem
implies that $s_{j,1}(b'),\ldots,s_{j,r_j}(b')$ are linearly
independent away from a thin subset.  

\begin{rem}
There are more precise versions of N\'eron's Theorem due to
Demyanenko, Manin and Silverman
(see \cite{silverman}, for example).
Masser has proposed another notion of what it means
for a subset of $B(K)$ to be small, 
known as `sparcity'.
For instance, the endomorphism ring of a family of
abelian varieties changes only
on a `sparse' set of rational points of the
base (see \cite{masser}).  For an analogue
to N\'eron's Theorem, see \cite{masser1}.
\end{rem}

Similar results hold for algebraic tori
and are proved using a version of N\'eron's Theorem for $\bG^n_m$-fibrations 
(see \cite{serre-MordellWeil} pp. 154).  
A sharper result (for 1-dimensional bases $B$)
can be obtained from the following recent theorem:

\begin{theo}
(\cite{bombieri-masser-zannier})
Let $C$ be an absolutely irreducible curve defined over
a number field  $K$ and $x_1,...,x_r$ rational functions in $
K(C)$, multiplicatively independent modulo constants.
Then the set of algebraic points $p\in C(\overline{\bQ})$
such that $x_1(p),...,x_r(p)$ are multiplicatively dependent
has bounded height.
\end{theo}

\

The main idea of the papers \cite{harris-tschi},
\cite{bog-tschi}, \cite{hassett-tschi} can be summarized as
follows.
We work over a number field $K$ and we assume that
all geometric data are defined over $K$.
Let $\pi : E\ra B$ be a Jacobian elliptic fibration
over a one dimensional base $B$.
This means that we have a family of curves of genus 1 and
a global zero section so that every fiber is in fact an
elliptic curve.
Suppose that we have another section $s$ which is
of infinite order in the Mordell-Weil group of
$E(K(B))$. The specialization results mentioned above show that
for a Zariski dense set of $b\in B(K)$ the restriction $s(b)$ is
of infinite order in the corresponding fiber $E_b$.
If $K$-rational points on $B$ are Zariski dense
then rational points on $E$ are Zariski dense as well.

Let us consider a situation when $E$ does not have any
sections but instead has a multisection $M$.
By definition, a multisection (resp. rational multisection)
$M$ is irreducible and the induced map $M\ra B$ is finite flat
(resp. generically finite) of degree $\deg(M)$.
The base-changed family $E\times_B M \ra M$ has
the identity section $\Id$ (i.e., the image of the diagonal under
$M\times_B M \ra E\times_B M$) and a (rational) section
$$\tau_M:=\deg(M)\Id - \Tr(M\times_B M)$$
where $\Tr(M\times_B M)$ is obtained (over the generic point)
by summing all the
points of $M\times_B M$.  By definition, $M$ is nondegenerate
if $\tau_M$ is nondegenerate.

When we are concerned only with rational points, we will ignore
the distinction between multisections and rational multisections,
as every rational multisection is a multisection over an
open subset of the base.  However, this distinction is crucial when
integral points are considered.

If $M$ is nondegenerate
and if rational points on $M$ are Zariski dense then rational
points on $E$ are Zariski dense
(see \cite{bog-tschi}).

\begin{exam}(\cite{harris-tschi})
Let $X$ be a quartic surface in $\bP^3$ containing a line $L$.
Consider planes $\bP^2$ passing through this line. The residual
curve has degree 3. Thus we obtain an elliptic fibration on
$X$ together with the trisection $L$. If $L$ is ramified
in a smooth fiber of this fibration then the multisection is
nondegenerate and rational points are
Zariski dense.
\end{exam}

This argument
generalizes to abelian fibrations $\pi : A\ra B$.
However, we do not know of any simple geometric
conditions insuring nondegeneracy of a (multi)section
in this case.
We do know that for any abelian variety
$A$ over $K$ there exists a finite extension
$K'/K$ with a nondegenerate point in
$A(K')$ (see \cite{hassett-tschi}).
This allows us to produce
nondegenerate sections over function fields.

\begin{prop}
Let $Y$ be a Fano threefold  of type $W_2$, that is
a double cover of $\bP^3$ ramified in a smooth
surface of degree 6.
Then rational points on the symmetric square $Y^{(2)}$
are potentially dense.
\end{prop}

{\em Proof.}
Observe that the symmetric square $Y^{(2)}$ is birational
to an abelian surface fibration over the Grassmannian of lines
in $\bP^3$. This fibration is visualized as follows:
consider two generic points in $Y$. Their images in $\bP^3$
determine a line, which intersects
the ramification locus in  6 points and
lifts to a (hyperelliptic) genus
two curve on $Y$. On $Y^{(2)}$ we have
an abelian surface fibration corresponding to
the degree 2 component of the relative Picard scheme.
Now we need to produce a nondegenerate multisection.
Pick two general points $b_1$ and $b_2$ on the branch surface.
The preimages in $Y$ of the corresponding tangent planes
are K3 surfaces $\Sigma_1$ and $\Sigma_2$,
of degree two with
ordinary double points at the points of tangency.
The surfaces $\Sigma_1$ and $\Sigma_2$ therefore have potentially
dense rational points (this was proved in \cite{bog-tschi}),
as does $\Sigma_1\times \Sigma_2$.  This is our
multisection;  we claim it is nondegenerate for
generic $b_1$ and $b_2$.
Indeed, it suffices to show that
given a (generic) point in $Y^{(2)}$, there exist $b_1$
and $b_2$ so that $\Sigma_1\times \Sigma_2$ contains
the point.  Observe that through a (generic) point
of ${\bP}^3$, there pass many tangent planes to the branch
surface.  \hfill   $\square$

\begin{rem}
Combining the above Proposition with the strong form
of Problem~\ref{prob:punc} we obtain
potential density of rational points on a Fano threefold
of type $W_2$ - the last family of smooth Fano threefolds
for which potential density is not known.
\end{rem}

\

Here is a formulation of the fibration method
useful for the analysis of integral points:
\begin{prop}\label{prop:group}
Let $B$ be a scheme over a number field $K$,
$G\ra B$ a flat group scheme,
$T\ra B$ an \'etale torsor 
for $G$, and $M\subset T$ a nondegenerate
multisection over $B$.  If $M$ has potentially
dense integral points
then $T$ has potentially dense integral points.
\end{prop}
{\em Proof.}
Without loss of generality, we may assume that
$B$ is geometrically connected and smooth.  The
base-changed family $T\times_B M$ dominates $T$, so
it suffices to prove density for $T\times_B M$.  Note
that since $M$ is finite and flat over $B$, $\tau_M$
is a well-defined {\em section} over all of $M$ (i.e., it is
not just a rational section).  Hence we may reduce
to the case of a group scheme $G\ra B$ with
a nondegenerate section $\tau$.

We may choose models $\cG$ and $\cB$ over $\Spec(\cO_S)$
so that $\cG\ra \cB$ is a group scheme with section $\tau$.
We may also assume that $\cO_S$-integral points of $\tau$
are Zariski dense.  The set of multiples $\tau^n$ of $\tau$, each
a section of $\cG\ra \cB$, is dense in $\cG$
by the nondegeneracy assumption.  
Since each has dense $\cO_S$-integral
points, it follows that $\cO_S$-integral points are Zariski
dense. \hfill  $\square$

A similar argument proves the following
\begin{prop}\label{prop:ellipK3}
Let $\varphi\,:\,X \ra \bP^1$ be a K3 surface with elliptic
fibration.  Let $M$ be a multisection over its image $\varphi(M)$,
nondegenerate and contained
in the smooth locus of $\varphi$.
Let $F_1,\ldots,F_n$ be fibers of $\varphi$ 
and $D$ a divisor supported
in these fibers and disjoint from $M$.  If $M$ has
potentially dense integral points then $(X,D)$ has potentially
dense integral points.
\end{prop}
{\em Proof.}
We emphasize that $X$ is automatically minimal
and the fibers of $\varphi$ are reduced
(see \cite{bog-tschi}).  Our assumptions imply that $M$ is
finite and flat over $\varphi(M)$.

After base-changing to $M$, we obtain a Jacobian elliptic fibration
$X':=X\times_{\bP^1}M$ with identity and a
nondegenerate section $\tau_M$.  Let $G\subset X'$ be the
open subset equal to the connected component of the identity.
Since $D':=D\times_{\bP^1}M$ is disjoint from the identity,
it is disjoint from $G$.  Hence it suffices to show that
$G$ has potentially dense integral points.

We assumed that $M$ is contained in the smooth locus of $\varphi$,
so $\tau_M$ is contained in the grouplike part of $X'$, and
some multiple of $\tau_M$ is contained in $G$.  Repeating
the argument for Proposition \ref{prop:group} gives the result.
\hfill $\square$

\section{Approximation techniques}
\label{sect:approx}

In this section we prove potential density of
integral points for certain pairs $(X,D)$ using
congruence conditions to control intersections
with the boundary.
Several of these examples are included
as support for the statement of
Problem~\ref{prob:punc}.

\begin{prop}\label{prop:group-punc}
Let $G=\prod_j^N G_j$
where $G_j$ are algebraic tori $\bG_m$ or
geometrically simple abelian varieties.
Let $Z$ be a subvariety in
$G$ of codimension
$>\mu=\max_j(\dim(G_j))$ and let $U=G\setminus Z$
be the complement. Then
integral points on $U$ are potentially
dense.
\end{prop}

{\em Proof.}  We are grateful to D. McKinnon
for inspiring the following argument.

The proof proceeds by induction on the number
of components $N$.  The base case $N=1$ follows
from the fact that rational points on
tori and abelian varieties are potentially dense,
so we proceed with the inductive step.  Consider
the projections $\pi'\,:\, G\ra G'=\prod_{j\neq N}G_j$
and $\pi_N\,:\,G \ra G_N$.
By assumption, generic fibers of $\pi'$ are
geometrically disjoint from $Z$.

Choose a ring of integers ${\cal O}_S$
and models ${\cal G}_j$ over $\Spec({\cal O}_S)$.
We assume that each ${\cal G}_j$ is smooth
over $\Spec({\cal O}_S)$ and that ${\cal G}_N$ has a
nondegenerate point
$q$ (see \cite{hassett-tschi}, for example,
for a proof of the existence
of such points on abelian varieties).

Let ${\cal T}$ be any subscheme of
${\cal G}_N$ supported over a finite subset of
$\Spec({\cal O}_S)$ such that ${\cal G}_N$
has an ${\cO}_S$-integral point $p_N$
disjoint from ${\cal T}$.
We claim that such integral points
are Zariski dense.
Indeed, for some $m>0$ we have
$$mq\equiv 0 \pmod{{\fp}}$$
for each ${\fp}\in \Spec(O_S)$ over which
${\cal T}$ has support.  Hence we may take the
translations of $p_N$ by multiples of $mq$.

After extending $\cO_S$, we may assume $U$ has at
least one integral point $p=(p',p_N)$ so that
${\pi'}^{-1}(p')$ and $\pi_N^{-1}(p_N)$ intersect
$Z$ in the expected dimensions.  In particular,
${\pi'}^{-1}(p')$ is disjoint from $Z$.
By the inductive hypothesis, we may extend
${\cO}_S$ so that
$$(\pi_N^{-1}(p_N)\simeq {\cal G}',
\pi_N^{-1}(p_N)\cap {\cal Z})$$
has dense integral points.  In particular,
almost all
such integral points are not contained in $\pi'({\cal Z})$,
a closed proper subscheme of ${\cal G'}$.
Let $r$ be such a point, so that
$F_r={\pi'}^{-1}(r)\simeq {\cal G}_N$ intersects
${\cal Z}$ in a subscheme ${\cal T}$ supported
over a finite number of primes.
Since $(r,p_N)\in F_r$ is disjoint
from ${\cal T}$, the previous claim
implies that the integral points of $F_r$
disjoint from ${\cal T}$ are Zariski dense.
As $r$ varies, we obtain a Zariski dense set of integral
points on ${\cal G}\setminus {\cal Z}$. \hfill $\square$

\begin{coro}
Let $X$ be a toric variety and $Z\subset X$
a subvariety of codimension $\ge 2$,
defined over a number field.
Then integral points on $(X,Z)$ are potentially dense.
\end{coro}

Another special case of the Arithmetic
puncturing problem~\ref{prob:punc} is
the following:

\begin{prob}
Are integral points on
punctured simple abelian varieties
of dimension $n>1$ potentially
dense?
\end{prob}

\begin{exam}
Potential density of integral points holds
for simple abelian varieties punctured in the origin,
provided that their ring of endomorphisms contains
units of infinite order.
\end{exam}

\section{Conic bundles and integral points}
\label{sect:conic}

Let $K$ be a number field, $S$ a finite set of places for $K$
(including all the infinite places), $\cO_S$ the corresponding
ring of $S$-integers, and $\eta\in \Spec(\cO_S)$ the generic point.
For each place $v$ of $K$, let $K_v$ be the corresponding
complete field and $\fo_v$ the discrete valuation ring
(if $v$ is nonarchimedean).  As before,
we use calligraphic letters (e.g., $\cX$)
for schemes (usually flat) over $\cO_S$ and roman letters (e.g., $X$)
for the fiber over $\eta$.

\subsection{Results on linear algebraic groups}
\label{subsection:gps}

Consider a linear algebraic group $G/K$.
Choose a model $\cG$
for $G$ over ${\cO}_S$, i.e., a flat group scheme of finite type
$\cG/{\cO}_S$ restricting to $G$ at the generic point.
This may be obtained by fixing a representation
$G \hookrightarrow \GL_n(K)$
(see also \cite{voskresenskii} \S 10-11).
The $S$-rank of $G$ (denoted $\rank(G,\cO_S)$) is defined as
the rank of the abelian group of sections of $\cG(\cO_S)$ over ${\cO}_S$.
This does not depend on the choice of a model.  
Indeed, consider two models
$\cG_1$ and $\cG_2$ with a 
birational map $b:\cG_1\dashrightarrow \cG_2$;
of course, $b$ is trivial over the generic point and the proper
transform of the identity section $I_1$ is the identity.
There is a subscheme $Z\subset \Spec(\cO_S)$ with
finite support such that the indeterminacy of $b$ is in the preimage
of $Z$.  It follows that the sections of $\cG_1$ congruent to $I_1$
modulo $Z$ have proper transforms which are sections of $\cG_2$.
Such sections form a finite-index subgroup of $\cG_1(\cO_S)$.

Let ${\bG}_m$ be the multiplicative group over ${\bZ}$,
i.e., $\Spec ({\bZ}[x,y]/\left<xy-1\right>)$;
it can be defined over an arbitrary scheme by extension of scalars.
There is a natural projection
$${\bG}_m({\bZ}) \ra \Spec({\bZ}[x])=\bA^1_{\bZ}\subset {\bP}^1_{\bZ}$$
so that ${\bP}^1_{\bZ} \setminus {\bG}_m(\bZ)=\{ 0,\infty \}$.
A {\em form} of ${\bG}_m$ over $K$
is a group scheme $G/K$ for which there
exists a finite field extension
$K'/K$ and an isomorphism $G\times_K K'\simeq {\bG}_m(K')$.
These are classified as follows (see \cite{Ono2} for a
complete account).  Any group automorphism
$$\alpha:{\bG}_m(K') \ra {\bG}_m(K')$$
is either inversion or the identity, depending on whether
it exchanges $0$ and $\infty$.
The corresponding automorphism group is smooth, so we may
work in the \'etale topology (see \cite{milne} Theorem 3.9).
In particular,
$$K-\text{forms of }{\bG}_m \simeq H^1_{\acute{e}t}(\Spec(K),{\bZ}/2{\bZ}).$$
Each such form admits a natural open imbedding into a projective
curve $G \hookrightarrow X$, generalizing the imbedding
of ${\bG}_m$ into ${\bP}^1$.  The complement $D=X\setminus G$
consists of two points.  The Galois action on $D$ is
given by the cocycle in
$H^1_{\acute{e}t}(\Spec(K),{\bZ}/2{\bZ})$ classifying $G$.

There is a general formula for the rank due
to T. Ono and J.M. Shyr (see \cite{Ono1}, 
Theorem 6 and \cite{Shyr}).
Let $T_v$ denote the completion of $T$
at some place $v$,
$\hat T$ and ${\hat T}_v$ the corresponding character groups,
and $\rho(T)$
(resp. $\rho(\cT_v)$) the number of independent elements
of ${\hat T}$ (resp. ${\hat \cT_v}$).  The formula takes the form
$${\rank}(T,\cO_S)=\sum_{v \in S}\rho(\cT_v)-\rho(T).$$
For forms of $\bG_m$ this is particularly simple.
For split forms
$$\rank(\bG_m,\cO_S)=\# \{ \text{places } v\in S\} - 1 .$$
Now let $G/K$ be a nonsplit form, corresponding 
to the quadratic extension
$K'/K$, and $S'$ the places of $K'$ lying over the places of $S$.
Then we have
$$\rank(G,\cO_S)=\# \{ \text{places } v \in S \text{ completely
splitting in } S' \}.$$

\subsection{Group actions and integral points}

Throughout this subsection, $\cX$ is a normal,
geometrically connected
scheme and $\cX \ra \Spec(\cO_S)$
a flat projective morphism.  Let $\cD\subset \cX$ be
an effective reduced Cartier divisor.  Contrary to our
previous conventions, we do not assume that $\cD$ is
flat over $\cO_S$.  Assume that a linear algebraic group
$G$ acts on $X$ so that $X\setminus D$ is a $G$-torsor.

\begin{prop}\label{prop:groupmodel}
There exists a model $\cG$ for $G$ such that
$\cG$ acts on $\cX$ and stabilizes $\cD$.
\end{prop}

{\em Proof.} Choose an imbedding $\cX \hookrightarrow
{\bP}^n_{\cO_S}$ and a compatible linearization
$G\hookrightarrow \GL_{n+1}(K)$ (see \cite{GIT}, Ch. 1 \S 3).
Let $\cG'\hookrightarrow \GL_{n+1}(\cO_S)$ be the
resulting integral model of $G$, so that $\cG'$ stabilizes the
ideal of $\cX$ and therefore acts on it.  Furthermore, $\cG'$
evidently stabilizes the irreducible components of $\cD$
dominating $\cO_S$.  The fibral components of $\cD$ are
supported over a finite subset of $\Spec(\cO_S)$.  We take
$\cG\subset\cG'$ to be the subgroup acting trivially over
this subset;  it has the desired properties.\hfill
$\square$

\begin{prop}\label{prop:infinitepts}
Assume $(\cX,\cD)$ has an $\cO_S$-integral point
and that $G$ has positive ${\cO}_S$-rank.
Then $(\cX,\cD)$ has an infinite number of ${\cO}_S$-integral
points.
\end{prop}

{\em Proof.}  Consider the action of ${\cG}(\cO_S)$ on
the integral point $\sigma$ (which has trivial stabilizer).
The orbit consists of
${\cO}_S$-integral points of $(\cX,\cD)$, an infinite
collection because $\cG$ has positive rank.
\hfill $\square$

Now assume that $X$ is a smooth rational curve.
A rational section (resp. bisection) $\cD \subset \cX$ is a
reduced effective Cartier divisor such that the generic fiber
$D$ is reduced of degree one (resp. two).
Note that the open curve $X\setminus D$ is
geometrically isomorphic to ${\bP}^1-\{ \infty \}$
(resp. ${\bP}^1-\{0,\infty\}$), and thus is a
torsor for some $K$-form $G$ of ${\bG}_a$ (resp. ${\bG}_m$).
This form is easily computed.  Of course,
${\bG}_a$ has no nontrivial forms.
In the ${\bG}_m$ case, we can regard $D_{\eta}$ as an element
of $H^1_{\acute{e}t}(\Spec(K),{\bZ}/2{\bZ})$, which gives the
descent data for $G$.

The following result is essentially due to Beukers
(see \cite{Beukers1}, Theorem 2.3):

\begin{prop}\label{prop:conicpoints}
Let $(\cX,\cD)\ra \Spec(\cO_S)$ be a rational curve
with rational bisection and $G$ the corresponding
form of ${\bG}_m$ (as described above).  Assume that
$(\cX,\cD)$ has an ${\cO}_S$-integral point and $\rank(G,{\cO}_S)>0$.
Then ${\cO}_S$-integral points of $(\cX,\cD)$ are Zariski dense.
\end{prop}

{\em Proof.}
This follows from Proposition~\ref{prop:infinitepts}.  Given
an ${\cO}_S$-integral point $\sigma$ of $(\cX,\cD)$, the
orbit $\cG(\cO_S)\sigma$ is infinite and thus Zariski dense.
\hfill  $\square$

Combining with the formula for the rank, we obtain the following:

\begin{coro}
Let $(\cX,\cD)\ra \Spec(\cO_S)$ be a rational curve
with rational bisection such that $(\cX,\cD)$
has an ${\cO}_S$-integral point.
Assume that either

a) $D$ is reducible over $\Spec(K)$ and $|S|>1$;  or

b) $D$ is irreducible over $\Spec(K)$ and at least
one place in $S$ splits completely in $K(D)$.

Then ${\cO}_S$-integral points of $(\cX,\cD)$ are Zariski dense.
\end{coro}

When $D$ is a rational section we obtain a similar result
(also essentially due to Beukers
\cite{Beukers1}, Theorem 2.1):
\begin{prop}\label{prop:linepoints}
Let $(\cX,\cD)\ra \Spec(\cO_S)$ be a rational curve
with rational section such that
$(\cX,\cD)$ has an ${\cO}_S$-integral point.
Then ${\cO}_S$-integral points of $(\cX,\cD)$ are Zariski dense.
\end{prop}

\subsection{$v$-adic geometry}
\label{subsection:vadic}

For each place $v\in S$,
consider the projective space ${\bP}^1(K_v)$ as a manifold
with respect to the topology induced by the $v$-adic
absolute value on $K_v$.  For simplicity, this will be
called the $v$-adic topology;  we will use the same term
for the induced subspace topology on ${\bP}^1(K)$.
Given an \'etale morphism of curves
$f:U \ra {\bP}^1$ defined over $K_v$,
we will say that $f(U(K_v))$
is a {\em basic \'etale open subset}.
These are open in the $v$-adic topology, either by
the open mapping theorem (in the archimedean case) or by Hensel's
lemma (in the nonarchimedean case).

Let
$$\chi_f(B):=\# \{ z \in \cO_{\{v\}}: |z|_v \le B \text{ and }
z \in f(U(K_v)) \}$$
where $B$ is a positive integer and
$$\cO_{\{v \}}:=\{ z \in K: |z|_w\le 1 \text{ for each }
w\ne v \}.$$
We would like to estimate the quantity
$$\mu_f:=\liminf_{B\ra \infty}\chi_f(B)/\chi_{\Id}(B)$$
i.e., the fraction of the integers contained in the image
of the $v$-adic points of $U$.

\begin{prop}\label{prop:etalecase}
Let $f:U\ra {\bP}^1$ be an \'etale morphism defined over $K_v$
and $f_1:C \ra {\bP}^1$ a finite morphism of smooth curves
extending $f$.  If there exists a point
$q\in f_1^{-1}(\infty)\cap C(K_v)$ at which $f_1$ is unramified
then $\mu_f=1$.
\end{prop}

{\em Proof.}  This follows from the fact that $f(U(K_v))$ is open
if $f$ is \'etale along $U$. \hfill $\square$

As an illustrative example, we take $K={\bQ}$ and $K_v=\bR$,
so that $\cO_{\{v\}}={\bZ}$.
The set $f(U(\bR))$ is a finite union of open intervals $(r,s)$ with $r,s\in
\bR\cup \{ \infty \}$, where the (finite) endpoints are
branch points.  We observe that
$$\mu_f=\begin{cases} 0 & \text{if $f(U(\bR))$ is bounded;}\\
    1/2& \text{if $\overline{f(U(\bR))}$ contains a
     one-sided neighborhood of $\infty$;}\\
       1& \text{if $\overline{f(U(\bR))}$ contains a
        two-sided neighborhood of $\infty$.}
        \end{cases}$$
We can read off easily which alternative occurs in terms of the
local behavior at infinity.  Let $f_1:C\ra {\bP^1}$
be a finite morphism of smooth curves extending $f$.
If $f_1^{-1}(\infty)$ has no real points
then $\mu_f=0$.
If $f_1^{-1}(\infty)$ has unramified (resp. ramified)
real points
then $\mu_f=1$ (resp. $\mu_f>0$.)

We specialize to the case of double covers:

\begin{prop}\label{prop:doublecase}
Let $U\ra {\bP}^1$ be an \'etale morphism defined over $K_v$
and $f_1:C \ra {\bP}^1$ a finite morphism of smooth curves
extending $f$.
Assume that $f_1$ has degree two and ramifies
at $q\in f_1^{-1}(\infty)$.  Then $\mu_f>0$.
\end{prop}

{\em Proof.}  Of course, $q$ is necessarily defined over $K_v$.
The archimedean case follows from the previous example,
so we restrict to the nonarchimedean case.
Assume $f_1$ is given by
$$y^2=c_nz^n+c_{n-1}z^{n-1}+\ldots+c_0,$$
where $z$ is a coordinate for the affine line in $\bP^1(K_v)$,
$c_n\ne 0$, and the $c_i\in  {\fo_v}$.
Substituting $z=1/t$ and
$y=x/t^{\lceil n/2 \rceil}$, we obtain the equation at
infinity
$$\begin{cases} x^2=c_n+c_{n-1}t+\ldots+c_0t^n & \text{for $n$ even}\\
x^2=c_nt+c_{n-1}t^2+\ldots+c_0t^n & \text{for $n$ odd}
\end{cases}.$$
If $n$ is even then $f_1^{-1}(\infty)$
consists of two non-ramified points, so we may
assume $n$ odd.  Then
$f_1^{-1}(\infty)$ consists of one ramification point
$q$, necessarily defined over $K_v$.

Write $c_n=u_0\pi^{\alpha}$ and $z=u_1\pi^{-\beta}$,
where $u_0$ and $u_1$ are units and $\pi$ is
a uniformizer in ${\fo_v}$.
(We may assume that some power $\pi^r$
is contained in $\cO_K$.)  Our equation takes the form
\begin{equation}
y^2\pi^{n\beta-\alpha}=u_0u_1^n+c_{n-1}u_1^{n-1}\pi^{\beta-\alpha}+
\ldots+c_0u_1 \pi^{n\beta-\alpha}\label{eqnsquare}.
\end{equation}

We review a property of the $v$-adic numbers, 
(proved in \cite{serre-LocalFields}, Ch. XIV \S 4).  
Consider the multiplicative group
$$U^{(m)}:=\{ u\in {\fo_v}: u\equiv 1\pmod{\pi^m} \}.$$
Then for $m$ sufficiently large we have
$U^{(m)}\subset K_v^2$.  In particular, to determine whether
a unit $u$ is a square, it suffices to consider its representative
$\mod{\pi^m}$.

Consequently, if $\beta$ is sufficiently large and has the same
parity as $\alpha$, then we can solve Equation \ref{eqnsquare}
for $y\in K_v$ precisely when $u_0u_1$ is a square.
For example, choose any $M\in \cO_K$ so that
$M\equiv u_0\pi^{(r-1)\beta}\pmod{\pi^{r\beta}}$ and set
$z=M/\pi^{r\beta} \in \cO_{\{ v \}}$.
Hence, of the
$z\in \cO_{\{v \}}$ with $|z|_v\le B$ (with $B\gg 0$), the fraction satisfying
our conditions is bounded from below.  It follows that $\mu_f>0$.
\hfill $\square$

Now let $f:U \ra {\bP}^1$ be an \'etale morphism of curves
defined over $K$.
Consider the function
$$
\omega_{f,S}(B):= \# \{ z \in \cO_S: |z|_v \le B
\text{ for each } v\in S \text{ and }
\alpha \in f(U(K)) \}
$$
and the quantity
$$\limsup_{B\ra \infty}\omega_{f,\{v\}}(B) / \chi_{f}(B).$$
We expect that this is zero provided that $f$ does not
admit a rational section.  We shall prove this is the case
when $f$ has degree two.

A key ingredient of our argument is a version of
Hilbert's Irreducibility Theorem:

\begin{prop}\label{prop:hilb}
Let $f:U \ra {\bP}^1$ be an \'etale morphism of curves, defined over $K$
and admitting no rational section.  Then we have
$$\limsup_{B\ra \infty}\omega_{f,\{v\}}(B) / \chi_{\Id}(B)=0.$$
\end{prop}

{\em Proof.}
We refer the reader
to Serre's discussion of Hilbert's irreducibility
theorem (\cite{serre-MordellWeil}, \S 9.6, 9.7).
Essentially the same argument applies in our situation.
\hfill $\square$

Combining Propositions
\ref{prop:etalecase}, \ref{prop:doublecase},
and \ref{prop:hilb}, we obtain:
\begin{coro}\label{coro:adicpoints}
Let $f:C \ra {\bP}^1$ be a finite morphism of smooth curves defined
over $K$.  Assume that $f$ admits no rational section and
that $f^{-1}(\infty)$ contains a $K_v$-rational point.
We also assume that $f$ has degree two.  Then we have
$$\limsup_{B\ra \infty}\omega_{f,\{v\}}(B) / \chi_{f}(B)=0.$$
In particular, the set
$\{ z \in \cO_{\{ v\}}: z \in f(C(K_v))\setminus f(C(K)) \}$
is infinite.
\end{coro}

\subsection{A density theorem for surfaces}

{\bf Geometric assumptions:}  Let $\cX$ and $\cB$ be flat and
projective over $\Spec(\cO_S)$ and
$\phi:\cX \rightarrow \cB$ be a morphism.  Let
$\cL\subset \cX$ be a closed irreducible
subscheme, $\cD \subset \cX$
a reduced effective Cartier divisor, and $\fq:=\cD\cap \cL$.
We assume the generic fibers satisfy the following:  $X$ is a
geometrically connected surface, $B$ a smooth curve,
$\phi:X\ra B$
a flat morphism such that the generic fiber is a rational
curve with bisection.
We also assume $L\simeq \bP^1_K$, $\phi|L$ is finite, and $L$
meets $D$ at a single point $q$, at which $D$ is nonsingular.
Write $\cX'$ for $\cX\times_{\cB}\cL$, $\cD'$ for $\cD\times_{\cB}\cL$,
$\cL'$ for the image of the diagonal in $\cX\times_{\cB}\cL$ (now
a section for $\phi':\cX'\rightarrow \cL$),
and $\fq'$ for $\cL'\cap \cD'$.
Finally, if $\cC'$ denotes the normalization of the union
of the irreducible components of $\cD'$ dominating $\cL$,
we assume that $\cC'\ra \cL$ has no rational section over $K$
(i.e., that $\cC'$ is irreducible over $K$).

\

\noindent {\bf Arithmetic assumptions:}
We assume that $(\cL,\fq)$
has an $\cO_S$-integral point.
Furthermore, we assume that for some $v \in S$, $C'$ has a $K_v$-rational point
lying over $\phi'(q')$.

\begin{rem}\label{rem:ramification}
This assumption is valid if any of
the following are satisfied:
\begin{enumerate}
\item{$D\ra B$ is unramified at $q$.}
\item{$D\ra B$ is finite (but perhaps ramified) at $q$
and $L\ra B$ has ramification
at $q$ of odd order.}
\item{$D\ra B$ is finite (but ramified) at $q$ and $L\ra B$ has
ramification at $q$ of order two.  Choose local uniformizers
$t,x,$ and $y$ so that we have local analytic equations
$t+ax^2=0$ and $t+by^2=0$  (with $a,b\in K$) for $D\ra B$
and $L\ra B$.
We assume that $ab$ is a square in $K_v$.}
\end{enumerate}
Note that in the last case, $D'$ and $C'$ have local analytic
equations $ax^2-by^2=0$ and $x/y=\pm\sqrt{b/a}$ respectively.
\end{rem}

\begin{theo}\label{theo:pelldense}
Under the geometric and arithmetic assumptions made above,
$\cO_S$-integral points of $(\cX,\cD)$ are Zariski dense.
\end{theo}

{\em Proof.}
It suffices to prove that
$\cO_S$-integral points of $(\cX',\cD')$ are Zariski dense.
These map to $\cO_S$-integral points $(\cX,\cD)$.

Consider first ${\cO}_S$-integral points of $(\cL',\fq')$.
These are dense by Proposition \ref{prop:linepoints},
and contain a finite index subgroup of $\bG_a(\cO_S)\subset
{\bP}^1_K$.  Corollary \ref{coro:adicpoints} and our
geometric assumptions imply
that infinitely many of these points lie in
$\phi'(C'(K_v)) \setminus \phi'(C'(K))$.

Choose a generic $\cO_S$-integral point $p$ of $(\cL',\fq')$
as described above.
Let $\cX'_p={\phi'}^{-1}(p),\cD'_p=\cX'_p\cap \cD',$ and
$\cL'_p=\cX'_p\cap \cL'$, so that $(\cX'_p,\cD'_p)$ is a rational
curve with rational bisection and integral point $\cL'_p$.
Combining the results of the previous paragraph with Proposition
\ref{prop:conicpoints}, with obtain that $\cO_S$-integral points of
$(\cX'_p,\cD'_p)$ are Zariski dense.  As we vary $p$, we
obtain a Zariski dense collection of integral points for
$(\cX',\cD')$.\hfill  $\square$

\subsection{Cubic surfaces containing a line}
Let $\cX_1$ be a cubic surface in ${\bP}^3_{\cO_S}$,
$\cD_1\subset \cX_1$
a hyperplane section, and $\cL_1\subset \cX_1$
a line not contained in $\cD_1$,
all assumed to be flat over $\Spec(\cO_S)$.
Write $\fq_1:=\cD_1 \cap \cL_1$, a rational section
over $\Spec(\cO_S)$.  Let
${\bP}^3_{\cO_S} \dashrightarrow \cB$ be the projection associated
with $\cL_1$, $\cX=\Bl_{\cL_1}\cX_1$, and
$\phi:\cX \ra \cB$ the induced projection
(of course, $\cB={\bP}^1_{\cO_S}$ if $\cO_S$ is a UFD).
Let $\cL\subset \cX$ be the proper transform of $\cL_1$,
$\cD\subset \cX$ the total transform of $\cD_1$,
and $\fq=\cL \cap \cD$.
We shall apply Theorem \ref{theo:pelldense} to obtain
density results for $\cO_S$-integral points of $(\cX_1,\cD_1)$.

We will need to assume the following geometric conditions:
\begin{enumerate}
\item[GA1]{$D_1$ is reduced everywhere and
nonsingular at $q_1$;}
\item[GA2]{$X_1$ has only rational double points as singularities,
with at most one singularity along $L_1$.}
\item[GA3]{$D_1$ is not the union of a line and
a conic containing $q_1$ (defined over $K$).}
\end{enumerate}
Using the first two assumptions, we analyze the projection
from the line $L_1$.
This induces a morphism
$$\phi:X \ra {\bP}^1.$$
Of course, $X=X_1$ if and only if $L_1$
is Cartier in $X_1$, which is the case exactly when
$X_1$ is smooth along $L_1$.  We use $L$ to denote
the proper transform of $L_1$
and $D$ to denote the proper transforms of
$L_1$ and $D_1$.
Our three assumptions imply that
$D$ equals the total transform of $D_1$ and
has a unique irreducible component $C$ dominating ${\bP}^1$.
We also have that
the generic fiber of $\phi$ is nonsingular,
intersects $D$ in two points, and
intersects $L$ in two points (if
$X_1$ is smooth along $L_1$) or in one point
(if $X_1$ has a singularity along $L_1$).
In particular, $L$ is a bisection (resp. section) of
$\phi$ if $X_1$ is nonsingular (resp. singular) along $L_1$.

We emphasize that
$\cO_S$-integral points
of $(\cX,\cD)$ map to $\cO_S$-integral points of $(\cX_1,\cD_1)$,
and all the Geometric Assumptions
of Theorem \ref{theo:pelldense} are satisfied except for the last one.
The last assumption is verified if any of the following
hold:
\begin{enumerate}
\item[GA4a]{The branch loci of $C\ra {\bP^1}$ and $L\ra {\bP}^1$
do not coincide.}
\item[GA4b]{The curve $C$ has genus one.}
\item[GA4c]{$X_1$ has a singularity along $L_1$.}
\end{enumerate}
Clearly, either the second or the third condition implies the first.

\

We turn next to the Arithmetic Assumptions.
\begin{enumerate}
\item[AA1]{
$(\cL_1,\fq_1)$ has an $\cO_S$-integral point.}
\end{enumerate}
Note that $\cO_S$-integral points of $(\cL_1,\fq_1)$ not lying
in the singular locus of $\cX_1\ra \Spec(\cO_S)$ lift naturally
to $\cO_S$-integral points of $(\cL,\fq)$.

Our next task is to translate the conditions of
Remark \ref{rem:ramification} to our situation.
They are satisfied in any of the following
contexts:
\begin{enumerate}
\item[AA2a]{$D_1$ is irreducible over $K$ and $q_1$ is not
a flex of $D_1$;}
\item[AA2b]{$X_1$ has a singularity along $L_1$;}
\item[AA2c]{$D_1$ is irreducible over $K$ and $q_1$ is
a flex of $D_1$.  Let $H$ be the hyperplane section containing
$L_1$ and the flex line.  We assume that $H\cap X_1=L_1\cup M$,
where $M$ is a smooth conic.}
\item[AA2d]{$D_1$ is irreducible over $K$ but $q$ is
a flex so that the hyperplane $H$ containing $L_1$ and the
flex line $F$ intersects $X_1$ in three coincident lines, i.e.,
$H\cap X_1= L_1\cup M_1\cup M_2$.  Choose local coordinates
$x$ and $y$ for $H$ so that $L_1=\{x=0 \}, F=\{y=0\},$
and $M_1\cup M_2= \{ ax^2+cxy+by^2=0 \}$.  Then we assume
that $ab$ is a square in $K_v$.}
\item[AA2e]{$D_1$ consists of a line and a conic $C_1$
irreducible over $K$, intersecting in two distinct points, each
defined over $K_v$.}
\end{enumerate}
In the first case, the map $D\ra B$ is unramified at $q$.
Note that in the second case $L$ is a section for $\phi$.
In the third case, our assumption implies that $L\ra B$
is unramified at $q$.  In the last case, we observe
that the points of $L$ lying over $\phi(q)$ are defined over $K$,
hence $C'$ has a $K_v$-rational point over $\phi'(q')$.

It remains to show that AA2d allows us to apply case 3 of
Remark \ref{rem:ramification}.  We fix projective coordinates
on ${\bP}^3$ compatibly with the coordinates already
chosen on $H$:
$y=0$ is the linear equation for the hyperplane containing
$D_1$, $z=0$ the equation for
$H$, $x=z=0$ the equations for $L_1$, and $x=z=w=0$ the equations
for $q_1$.  Under our assumptions,
the equations for $D_1$ and $X_1$ take the form
\begin{eqnarray*}
g&:=&zw^2+ax^3+c_1wxz+c_2wz^2+c_4x^2z+c_5xz^2+c_6z^3=0\\
f&:=&g+cx^2y+bxy^2+yz\ell(w,x,y,z)=0
\end{eqnarray*}
where $\ell$ is linear in the variables.
The conic bundle structure $\phi:X \ra B$ is obtained
by making the substitution $z=tx$
\begin{eqnarray*}
g'&=&tw^2+x(wc_1t+wc_2t^2)+x^2(a+c_4t+c_5t^2+c_6t^3)=0\\
f'&=&g'+cxy+by^2+ty\ell(w,x,y,tx).
\end{eqnarray*}
We analyze the local behavior of $D\ra B$ at $q$
using $x$ as a coordinate for $D$.
First dehomogenize
$$g''=t+x(c_1t+c_2t^2)+x^2(a+c_4t+c_5t^2+c_6t^3)=0$$
and then take a suitable analytic change of coordinate on $D$
to obtain $t+aX^2=0$.  To analyze $L\ra B$, we set $x=0$
and use $y$ as a coordinate
$$f''=t+by^2+ty\ell(1,0,y,0)=0.$$
After a suitable analytic change of coordinate on $L$, we
obtain $t+bY^2=0$.

\begin{rem}
\label{rem:realcase}
We further analyze condition AA2d when
$K_v={\bR}$.  Then $ab$ is a square if and only
if $ab\ge 0$.  This is necessarily the case if $c^2-4ab<0$,
i.e., if the lines $M_1$ and $M_2$ are defined over an
imaginary quadratic extension.
\end{rem}

We summarize our discussion in the following theorem:
\begin{theo}\label{theo:cubicdense}
Let $\cX_1$ be a cubic surface,
$\cD_1\subset \cX_1$
a hyperplane section, and $\cL_1\subset \cX_1$
a line not contained in $\cD_1$,
all assumed to be flat over $\Spec(\cO_S)$.
Write $\fq_1:=\cD_1 \cap \cL_1$.
Assume the following:
\begin{enumerate}
\item{GA1,GA2,GA3, and AA1;}
\item{at least one of the assumptions GA4a,GA4b,or GA4c;}
\item{at least one of the assumptions AA2a,AA2b,AA2c,AA2d, or AA2e.}
\end{enumerate}
Then $\cO_S$-integral points of $(\cX_1,\cD_1)$ are Zariski dense.
\end{theo}

\

We recover the following result (essentially Theorem 2 of Beukers
\cite{Beukers2}):
\begin{coro}\label{coro:BeukersThm2}
Let $\cX_1$ be a cubic surface,
$\cD_1\subset \cX_1$
a hyperplane section, and $\cL_1\subset \cX_1$
a line not contained in $\cD_1$,
all assumed to be flat over $\Spec(\bZ)$.
Write $\fq_1:=\cD_1 \cap \cL_1$.
Assume that
\begin{enumerate}
\item{$X_1$ and $D_1$ are smooth;}
\item{there exists an ${\bZ}$-integral point of $(\cL_1,\fq_1)$;}
\item{if $q$ is a flex of $D_1$, we assume that the hyperplane
containing $L_1$ and the flex line intersects $X_1$ in
a smooth conic and $L_1$.}
\end{enumerate}
Then $\bZ$-integral points of $(\cX_1,\cD_1)$ are Zariski dense.
\end{coro}

We also recover a weak
version of Theorem 1 of \cite{Beukers2}.
(This theorem is asserted to be true but the proof is not
quite complete;  the problem occurs in the argument
for the second part of Lemma 2.)
\begin{coro}\label{coro:BeukersThm1}
Retain all the hypotheses of Corollary \ref{coro:BeukersThm2},
except that we allow the existence of a hyperplane $H$
intersecting $X_1$ in three lines $L_1,M_1,$ and $M_2$
and containing a flex line $F$ for $D_1$ at $q$.
Let $p$ be a place for $\bZ$ (either infinite or finite).
Choose local coordinates
$x$ and $y$ for $H$ so that $L_1=\{x=0 \}, F=\{y=0\},$
and $M_1\cup M_2= \{ ax^2+cxy+by^2=0 \}$, and assume
that $ab$ is a square in ${\bQ}_p$.  Then $\bZ[1/p]$-integral points
of $(\cX_1,\cD_1)$ are Zariski dense (where $\bZ[1/\infty]={\bZ}$
and ${\bQ}_{\infty}=\bR$.)
\end{coro}
Of course, there are infinitely many primes $p$ such that $ab$ is a
square in $\bQ_p$.  When $p=\infty$,
by Remark \ref{rem:realcase}
it suffices to verify that $M_1$ and $M_2$ are defined over
an imaginary quadratic extension.

We also obtain results in cases where the boundary is
reducible:

\begin{coro}
Let $\cX_1$ be a cubic surface,
$\cD_1\subset \cX_1$
a hyperplane section, and $\cL_1\subset \cX_1$
a line not contained in $\cD_1$,
all assumed to be flat over $\Spec(\bZ)$.
Write $\fq_1:=\cD_1 \cap \cL_1$.
Assume that
\begin{enumerate}
\item{$X_1$ is smooth;}
\item{there exists an $\cO_S$-integral point of $(\cL_1,\fq_1)$;}
\item{$D_1=E \cup C$, where $E$ is a line intersecting $L_1$ and
$C$ is a conic irreducible over $K$;}
\item{$C$ intersects $E$ in two points, defined over $K_v$
where $v$ is some place in $S$;}
\item{there exists at most one conic in $X_1$ tangent
to both $L_1$ and $C$.}
\end{enumerate}
Then $\cO_S$-integral points of $(\cX,\cD)$ are Zariski dense.
\end{coro}

Note that the assumption on the conics tangent to $L_1$ and $C$
is used to verify GA4a.

\subsection{Other applications}

Theorem \ref{theo:pelldense} can be applied in many
situations.  We give one further example:

\begin{theo}\label{theo:P1P1}
Let $\cX={\bP}^1_{\cO_S} \times {\bP}^1_{\cO_S}$,
$\cD\subset \cX$ a divisor of type $(2,2)$, and $\cL\subset \cX$
a ruling of $\cX$, all flat over $\cO_S$.
Assume that
\begin{enumerate}
\item{$D$ is nonsingular;}
\item{$L$ is tangent to $D$ at $q$;}
\item{$\cO_S$-integral points of $(\cL,\fq)$ are
Zariski dense.}
\end{enumerate}
Then $\cO_S$-integral points of $(\cX,\cD)$ are Zariski
dense.
\end{theo}

{\em Proof.} Let $\phi$ be the projection for which $\cL$
is a section.  Since $\cC=\cD$ in this case, the second arithmetic
assumption of Theorem \ref{theo:pelldense} is easily satisfied.
\hfill $\square$

\section{Potential density for log K3 surfaces}

We consider the following general situation:

\begin{prob}[Integral points of
log K3 surfaces]\label{prob:xd}
Let $X$ be a surface and $D$ a reduced effective
Weil divisor such that $(X,D)$ has log terminal singularities
and $K_X+D$ is trivial.
Are integral points on $(X,D)$ potentially dense?
\end{prob}

Problem \ref{prob:xd} has been studied when $D=\emptyset$
(see, for example, \cite{bog-tschi}).
In this case density holds if $X$ has infinite automorphisms
or an elliptic fibration.

The case $X={\bP}^2$ and $D$ a plane cubic has also attracted
significant attention.  Silverman \cite{silverman-2} proved
potential density in the case where $D$ is singular and raised
the general case as an open question.  Beukers \cite{Beukers1}
established this by considering the cubic surface $X_1$ obtained
as the triple cover of $X$ totally branched over $D$.

Implicit in \cite{Beukers2} is a proof of potential density when
$X_1$ is a smooth cubic surface and
$D_1$ is a smooth hyperplane section.
Note that this also follows from Theorem \ref{theo:cubicdense}
(cf. also Corollaries
\ref{coro:BeukersThm2} and \ref{coro:BeukersThm1}.)
After suitable extensions of $K$ and additions to $S$,
there exists a line $L\subset X$ defined over $K$ and
the relevant arithmetic assumptions are satisfied.
Similarly, the case of $X={\bP}^1\times {\bP}^1$ and $D$ a smooth
divisor of type $(2,2)$ follows from Theorem \ref{theo:P1P1}.

More generally, let $X$ be a smooth Del Pezzo surface of index one,
i.e., $K_X{\bZ}$ is saturated in $\Pic(X)$, and
with degree $d:=K_X^2\ge 4$.  Let $D$ be a smooth anticanonical
divisor.  Choose general points
$W=\{x_1,\ldots,x_{d-3} \} \subset X\setminus D$, and let
$X_1=\Bl_W X$ and $D_1$ be the proper transform of $D$.
Hence $X_1$ is a cubic surface, $D_1$
is a smooth hyperplane section, and the
induced map of pairs
$$(X_1,D_1) \ra (X,D)$$
is dominant.  Since integral points for $(X_1,D_1)$ are
potentially dense, the same holds true for $(X,D)$.

We summarize our results as follows:
\begin{theo}\label{theo:logK3dense}
Let $X$ be a smooth Del Pezzo surface of degree $\ge 3$
and $D$ a smooth anticanonical divisor.  Then integral
points for $(X,D)$ are potentially dense.
\end{theo}

\

We close this section with a list of open special cases of
Problem \ref{prob:xd}.

\begin{enumerate}
\item{Let $X$ be a Del Pezzo surface of degree one or two
and $D$ an anticanonical cycle.  Show that integral points
for $(X,D)$ potentially dense.}
\item{Let $X$ be a Hirzebruch surface and $D$
an anticanonical cycle.  Find a smooth rational
curve $L$, intersecting $D$ in exactly one point $p$,
so that the induced map $\varphi:L\ra {\bP}^1$
is finite surjective.}
\end{enumerate}

\subsection{Appendix: some geometric remarks}

The reader will observe that the methods employed to
prove density for integral points on conic bundles
(with bisection removed) are not quite analogous to the
methods used for elliptic fibrations. The 
discrepancy can be seen in a number of ways. First, given a
multisection $M$ for a conic bundle (with bisection removed),
we can pull-back the conic bundle to the multisection.
The resulting fibration has two {\em rational} sections,
$\Id$ and $\tau_M$ (see section \ref{sect:fibration}).
However, {\em a priori} one cannot control how
$\tau_M$ intersects the boundary divisor (clearly, this is
irrelevant if the boundary is empty).
A second explanation may be found in the lack of a good
theory of (finite type) N\'eron models for algebraic
tori (see chapter 10 of \cite{blr}).

We should remark that in some special cases these difficulties
can be overcome, so that integral points may be obtained
by geometric methods completely analogous to those 
used for rational points.  Consider the cubic surface
$$x^3+y^3+z^3=1$$
with distinguished hyperplane at infinity.  This surface
contains a line with equations $x+y=z-1=0$.  Euler showed
that the resulting conic bundle admits a multisection
$(x_0,y_0,z_0)=(9t^4,3t-9t^4,1-9t^3),$
which may be reparametrized as
$(x_1,y_1,z_1)=(9t^4,-3t-9t^4,1+9t^3).$
Lehmer \cite{lehmer} showed that this is the first in a sequence of
multisections, given recursively by
\begin{eqnarray*}
(x_{n+1},y_{n+1},z_{n+1})&=&
2(216t^6-1)(x_n,y_n,z_n)-(x_{n-1},y_{n-1},z_{n-1}) \\
&+&(-108t^4,-108t^4,216t^4+4)
\end{eqnarray*}
This should be related to the fact
that the norm group scheme
$$u^2-3(108t^6-1)v^2=1,$$
admits a section of infinite order $(u,v)=(216t^6-1,12t^3)$.

\end{document}